\documentclass[12pt]{article}
\usepackage{amsfonts,amssymb,amsxtra}
\usepackage{url}
\usepackage{amsmath,amscd}
\usepackage{ntheorem}
\def\eps{{\boldsymbol\epsilon}}

\newtheorem*{proposition-non}{Proposition}
\newtheorem{theorem}{Theorem}[section]
\newtheorem{corollary}[theorem]{Corollary}
\newtheorem{conjecture}[theorem]{Conjecture}
\newtheorem{lemma}[theorem]{Lemma}
\newtheorem{proposition}[theorem]{Proposition}

\newtheorem{remark}[theorem]{Remark}
\newtheorem{example}[theorem]{Example}
\newtheorem{definition}[theorem]{Definition}
\newtheorem{notation}[theorem]{Notation}
\DeclareSymbolFont{AMSa}{U}{msa}{m}{n}
\DeclareMathSymbol{\boxtimes}{\mathbin}{AMSa}{"02}
\def\fin{\nolinebreak $\Box$}
\def\proof{{\noindent\it Proof: }}

\def\dual{ ^{\vee}}
\def\LL{{\cal L}}

\def\quot#1#2{\mbox{\raise 0.1cm\hbox{$#1$}{\big /}\raise -0.1cm\hbox{$#2$}}}

\def\rela{{\rm Z} \kern -0.17cm \raise 0.05cm\hbox{\tiny{\it /}}\kern
  0.05cm}
\def\petitrela{{\rm Z} \kern -0.11cm \raise 0.05cm\hbox{{\mini /}}}
\def\ratio{{\rm Q} \kern -0.216cm\rule{0.25pt}{8pt}\kern 0.216cm}
\def\petitratio{{\rm Q} \kern -0.16cm\rule{0.15pt}{5pt}\kern 0.13cm}
\def\petitcomplex{{\rm C} \kern -0.13cm\rule{0.15pt}{5pt}\kern 0.13cm}
\def\PP{\mathbb{P}}
\def\complex{\mathbb{C}}
\newcommand{\OO}{{\cal O}}

\newlength{\elargir}
\title{Pfaffian bundles on cubic surfaces and configurations of planes}
\author{Han Fr\'ed\'eric}

\begin{document}

\maketitle

 %\phantom{.}
 %\vskip1cm 
 %\noindent{\bf \LARGE }
% \vskip1cm
 %\noindent{\bf Han Fr\'ed\'eric}
 \indent  Institut de Math\'ematiques de Jussieu - Paris Rive Gauche,
 
           Universit\'e Paris 7 - 5 rue Thomas Mann,

           Batiment Sophie-Germain

           75205 Paris Cedex 13, FRANCE

           email: han@math.jussieu.fr


\bigskip
% %
\begin{abstract}
We give a canonical birational map between the moduli space of pfaffian vector bundles on
a cubic surface and the space of complete pentahedra inscribed in the cubic
surface. The universal situation is also considered, and we obtain a
rationality result. As a by-product, we provide an explicit normal form for five general lines in
$\PP_5$. Applications to the geometry of  Palatini threefolds and 
Debarre-Voisin's Hyper-K\"ahler manifolds are also discussed. 
\end{abstract}

\def\P5{{\mathbb{P}_5}}
\def\Pw{\mathop{\pi_3}}
\def\Pwv{\mathop{\pi_3 ^\vee}}

\def\Proj#1{\mathop{Proj(S^\bullet(#1))}}
\def\t#1{\mathop{^t\! #1}}

\section{Introduction}
In the classical theory of determinantal hypersurfaces, the case of pfaffian
cubics has already found many applications. (\cite{Ma-Ti}, \cite{I-R},
\cite{Dr}). This paper presents new invariants for these objects.

Let $\P5$ be a five dimensional projective space over the complex numbers, and
denote by $V_6$ the vector space $H^0 (\OO_{\P5}(1))$. For $n\geq 2$, let
$\pi_{n-1}$ be a projective space of dimension $n-1$, and
$W_n=H^0(\OO_{\pi_{n-1}}(1))$.
\begin{definition}\label{defpfaf}
  For $2\leq n$, a general element of $\bigwedge\limits^2 V_6 \otimes W_n$ gives
  a skew-symmetric map $M$ with linear coefficients over $\pi_{n-1}$. We have the exact sequence:  
   $$\begin{CD}
   0@>>>  V_6^\vee \otimes \OO_{\pi_{n-1}}(-1) @>M>> V_6 \otimes \OO_{\pi_{n-1}} @>>> F@>>> 0 
 \end{CD}$$
%\shortE{V_6^\vee \otimes \OO_{\P5}(-1)}{V_6 \otimes \OO_{\P5}}{F}{M}
where $F$ is a rank $2$ sheaf over the pfaffian divisor. For $n\leq 6$, the sheaf $F$  is a vector
bundle over a smooth cubic. Let's  call it the pfaffian bundle defined by $M$.
\end{definition}

It is known from  classical works on  representations of a cubic form by
pfaffians (Cf \cite{Be}, \cite{Do}),  that for $n\geq 6$ a general cubic
divisor is not a pfaffian, and for $3\leq n \leq 5$ the pfaffian bundles have
moduli spaces of strictly positive dimension.  

\bigskip
The main result of this article concerns the $n=4$ case. In this situation we
have the following results from \cite{Be}:
\begin{itemize}
\item[-] Every smooth cubic surface of $\Pw$ can be defined by a linear pfaffian.
\item[-] Let $(W_4\otimes \bigwedge\limits^2 V_6)^{sm}$  be the openset of
  $W_4\otimes \bigwedge\limits^2 V_6$ corresponding to  smooth pfaffian
  surfaces. For any element $M$ of   $(W_4\otimes \bigwedge\limits^2 V_6)^{sm}$,
  the pfaffian sheaf is a stable rank $2$ vector bundle over the pfaffian surface
  $Pf(M)$. Moreover, it is an arithmetically Cohen-Macaulay
  sheaf, and every arithmetically Cohen-Macaulay rank $2$ vector
  bundle over a smooth cubic surface $S$ with determinant $\OO_S(2)$ is a
  pfaffian bundle.
\item[-]  The quotient of $(W_4\otimes
  \bigwedge\limits^2 V_6)^{sm}$ by $GL(V_6)$ with the following action:  
  $$\begin{array}[c]{ccc}
    GL(V_6)\times (W_4\otimes \bigwedge\limits^2 V_6)^{sm} & \to & (W_4\otimes
    \bigwedge\limits^2 V_6)^{sm} \\
    (P,M) & \mapsto &  ^t P.M.P
  \end{array}.$$
is isomorphic to the space of pairs $(S,F)$ where $S$ is a smooth cubic surface
of $\Pw$ and $F$ an isomorphism class of a pfaffian bundle on $S$. It is a geometric quotient.
\end{itemize}
In this article we obtain a geometric interpretation of these orbits.
%Briefly, we have the following: 
\begin{definition}\label{inscribed}
A complete pentahedron inscribed in a cubic surface $S$ of the projective space
$\Pw $ is a set $\{H_0, \dots, H_4\}$ of $5$ planes  of $\Pw $ such that:
\begin{itemize}
\item[i)] $(H_0, \dots, H_4)$ is a projective basis of $\Pwv$.
\item[ii)] The $10$ points $(H_i\cap H_j \cap H_k)_{0\leq i<j<k\leq 4}$ are on $S$. 
\end{itemize}
We define the subset $\mathcal{H}$ of $|O_{\Pw}(3)|\times |O_{\Pw}(5)|$
(resp. the subset $\mathcal{H}_{ord}$ of $|O_{\Pw}(3)|\times (\Pw)^5$) to be the set of elements
$(S,\Pi)$ such that $S$ is a smooth cubic surface of $\Pw$ and $\Pi$ is a complete
pentahedron inscribed in $S$ (resp. complete pentahedron inscribed in $S$ with an
order on the five planes). For a cubic surface $S$ of $\Pw$, denote by
$\mathcal{H}_S$ the pullback in $\mathcal{H}$ of $S$ by the first projection.
\end{definition}

The first four sections give two natural methods to construct $5$ hyperplane
sections of a cubic surface from a pfaffian vector bundle. Eventually they are
generically identical, and we obtain:
\newpage
\begin{theorem}\label{mainbir}
  There is a natural birational map from $(W_4\otimes \bigwedge\limits^2
  V_6)^{sm}/GL(V_6)$ to $\mathcal{H}$  such that the composition with the
  projection to $|\OO_{\Pw}(3)|$ is the pfaffian map. In particular:
  \begin{itemize}
  \item[.] $(W_4\otimes \bigwedge\limits^2 V_6)^{sm}/GL(V_6)$ is a rational
    variety of  dimension $24$.
  \item[.] Let  $S$ be a general cubic surface. The moduli space of pfaffian
    bundles on $S$ is birational to $\mathcal{H}_S$.  
  \end{itemize}
\end{theorem}
%Additionnally, the first four sections enlight this statement with two natural
%ways to construct $5$ hyperplane sections of a cubic surface from a
%pfaffian vector bundle. 
Both constructions enlight this theorem differently. The first one: $\Phi_1$
(cf definition \ref{defPHI1})   is a classical 
problem of hyperplane restriction of the pfaffian bundles. So section 
\ref{secn3} starts with the easy case $n=3$ to introduce some invariants
of these bundles. The universal situation is then described because
many geometric interpretations of the later sections are specializations 
of this construction.

\bigskip
Section \ref{secn4} details the $n=4$ case to settle the construction of $\Phi_1$. The
projectivisation of a pfaffian bundle on a cubic surface is called a Palatini threefold.
Such varieties are well-known to be the only known examples of
smooth $3$-dimensional varieties $X$ in $\P5$ such that
$h^0(\OO_X(2))>h^0(\OO_{\P5}(2))$. One can find references for their Hilbert scheme
(\cite{Fa-Fa}, \cite{Fa-Me}), and also  references where they are
in a list of exceptions to some geometrical property (cf \cite{Me-Po},
\cite{Ot}). Some of their classical properties are also described in \cite{Do}
and \cite{Ok}. But here some new results are required. First we give an interpretation of
their anticanonical linear system to prove that it is $\Pwv$. Then we describe
this linear system in proposition \ref{propantican}. It turns out that its exceptional
locus is $5$ points of $\Pwv$. This achieves the construction of $\Phi_1$.

\bigskip
But the geometric configuration of these planes is only explained in section \ref{W2V6} by
the second construction: $\Phi_2$ (cf corollary \ref{defPHI2}).  This time, it
is a problem of sum of matrices of rank $2$. The key step to construct $\Phi_2$
is the surprising proposition \ref{P4-5sec} with  following summary: 
\begin{proposition-non}
The projection of the Grassmannian of lines $G(2,V_6)$ from a general
$3$-dimensional projective space has a single point of order $5$.   
\end{proposition-non}
The claim that $\Phi_1$ and $\Phi_2$ are generically the same, and also
their birationality, are  proved in section
\ref{sectionexplicit}  from the explicit formula of theorem\ref{explicit}. This
ends the proof of theorem \ref{mainbir}. As a 
by-product we obtain in corollary \ref{unirationalisation} an explicit
generically finite unirationalization of the quotient of the product of five copy of $G(2,V_6)$
by the diagonal action of $PGL(V_6)$. 

\smallskip
Recently, 
F. Tanturri (cf \cite{T}) found an algorithm to obtain a pfaffian representation from the
equation of a cubic surface. Although some representations are similar, the main
difference is that any pfaffian bundle on the surface would solve his problem, while in
our situation we have additional requirements such that only one bundle is
solution. 

%%%%%
%%%%
\bigskip
In the last section we investigate those properties over a base. We explain how
the Debarre-Voisin's symplectic manifold can be considered as a parameter space
for Palatini threefolds in a  six dimensional variety of $\PP_9$. Those 
varieties of dimension six were discovered by C. Peskine. They  are of independent interest
because they are smooth and non quadratically normal in $\PP_9$ (it's a boundary
case in Zak's theory of quadratic normality). However, most
of their geometric properties are unknown. In particular, it would be very
interesting to understand those varieties from a Palatini threefold in a similar
way that a Veronese surface is related to $\PP_2\times \PP_2$. So we will also
explain in this section the consequences on the Peskine's varieties of the work
on the Palatini threefolds done in  section \ref{secn4}.

\begin{description}
\item[\bf Aknowledgement:] \ \\
  I'd like to thanks I. Dolgachev for encouraging discussions and references.
\end{description}

\section{Invariants of Pfaffian bundles over plane cubics.}\label{secn3}
\subsection{Ruled surfaces in $\P5$, and the $n=3$ case.}
In this section, we detail the case $n=3$. The following easy lemma is a basic
step that enlights the next sections.
\begin{lemma}\label{planecubi}
  For a general element of $W_3\otimes \bigwedge\limits^2 V_6$, we consider the
  associated exact sequence:  
  \begin{equation}\label{pfaffP2}
   \begin{CD}
    0@>>>  V_6^\vee \otimes \OO_{\pi_{2}}(-1) @>M>> V_6 \otimes \OO_{\pi_{2}} @>>> F@>>> 0 
   \end{CD}
  \end{equation}

%\shortE{V_6^\vee \otimes \OO_{\P5}(-1)}{V_6 \otimes \OO_{\P5}}{F}{M}
with $M= -\t{M}$. The cokernel  $F$ is a rank $2$ vector bundle over the smooth plane
cubic $C$ defined by the pfaffian of $M$, and $F$ is isomorphic to one of the following
bundles: 
\begin{itemize}
\item[a)] $\LL(1) \oplus \LL\dual(1)$, where $\LL$ is a line bundle of degree $0$ on $C$ such
  that $h^0(\LL^2)=0$.
\item[b)] $F$ is the unique unsplit extension: $$0 \rightarrow \theta(1) \rightarrow
  F\rightarrow \theta(1)\rightarrow 0 $$ where $\theta^2= \OO_C$ and $\theta\neq \OO_C$.
\item[c)] $F=\theta(1) \oplus \theta(1)$ where $\theta^2= \OO_C$ and $\theta\neq \OO_C$.
\end{itemize}
\end{lemma}
\proof
To simplify the notations, let $F_0$ denote $F(-1)$.
First one can remark that $h^0(F_0)=0$, and that $F_0 \simeq (F_0)\dual$ because $M$
is skew-symmetric. So we have $\wedge^2(F_0) = \OO_C $. We  choose a point $p$ on $C$. We
will now  prove that there  is 
a point $r$ of $C$ such that $h^0(F_0(p-r))>0$.  

From Riemann-Roch's theorem the bundle $F_0(p)$ has a pencil
of sections. This gives,  on $\PP_1\times C$, a section of the bundle $\OO_{\PP_1}(1)
\boxtimes F_0(p)$. But the computation of the second Chern's class of this bundle implies
that this section has a non empty vanishing locus, so there is a point $r$ of $C$ such
that $h^0(F_0(p-r))>0$. Let's recall that  $h^0(F_0)=0$ to obtain that $\OO_C(p-r)$ is
not trivial  and that $F$ is isomorphic to one of the $3$  above cases. \fin

\begin{remark}\label{3casplans}
The ruled surface $\Proj{F}$ has a natural embedding in $\P5$ given by the surjection
  in the sequence (\ref{pfaffP2}) such that in the cases:
\begin{itemize}
\item[a)] it contains $2$ plane cubic curves, and the  planes spanned by these curves
  are disjoint in $\P5$.
\item[b)] it contains only one plane cubic.
\item[c)] it contains infinitely many plane cubics.  The planes spanned by these curves
  are the planes of a Segre: $\PP_1\times \PP_2 \subset \P5$.
\end{itemize}
Moreover, the planes in those $3$ cases are the planes of $\P5$ isotropic for all the
skew-symmetric forms defined by $M$.  
\end{remark}
\proof In those $3$ cases, the bundle $F$ has an invertible quotient of rank $1$ and degree
$3$. We just have to show that those embeddings of $C$ are isotropic for $M$. But it is a
corollary of the fact that the resolution of $F$ can have a skew-symmetric form deduced from
the isomorphism: $\wedge^2(F(-1)) \simeq \OO_C$. Conversely, any isotropic plane for $M$ gives
the existence of $P\in GL(V_6)$ such that: $\t{P}.M.P=
\left(\begin{array}{cc}
  0 & -\t{A}\\
  A & B
\end{array}\right)
$, where $A,B$ are $3$ by $3$ matrices with linear entries. So the cokernel of $A$ gives the expected invertible quotient of $F$ of degree $3$. \fin

\subsection{Universal settings and the $SL(V_6)$-invariant double cover}
\begin{definition}
Let $G(3,V_6\dual)$ and $G(3,\bigwedge\limits^2 V_6)$  be the Grassmannians of
$3$-dimensional vector subspaces of $V_6\dual$ and  $\bigwedge\limits^2 V_6$. Denote
by $K_3$ and $R_3$ their tautological subbundles. We define the isotropic
incidence:
$$
\begin{CD}
Z\subset G(3,V_6\dual) \times G(3,\bigwedge\limits^2 V_6)  @>p_2>>
G(3,\bigwedge\limits^2 V_6)\\
@Vp_1VV \\
G(3,V_6\dual)
\end{CD}
$$
to be the vanishing locus of the unique $SL(V_6)$-invariant section of $\bigwedge\limits^2
K_3\dual \boxtimes R_3\dual$. Denote by $\mathcal{U}$ the open subset of $G(3,
\bigwedge\limits^2 V_6)$ made of  subspaces  such that the intersection of
their projectivisation with the pfaffian hypersurface of $\PP(\bigwedge\limits^2 V_6)$ is a
smooth cubic curve.

The restriction of $Z$ to $G(3,V_6\dual)\times \mathcal{U}$ will be noted: $Z_{\mathcal{U}}$. Let $E_{12}$ be
the rank $12$ bundle defined by the exact sequence:

\begin{equation}
  \label{eqE12}
\begin{CD}
   0@>>>  E_{12} @>>> \bigwedge\limits^2 V_6 \otimes \OO_{G(3,V_6\dual)} @>>> \bigwedge\limits^2 K_3\dual @>>> 0 
 \end{CD}
\end{equation}
\end{definition}

I'd like to thanks A. Kuznetsov for the following description of $Z$ from the relative Grassmannian. 
\begin{proposition}\label{propisoE12}
The isotropic incidence $Z$ is isomorphic to the relative Grassmannian 
$G(3,E_{12})$ of linear subspaces of the bundle $E_{12}$. The projection $Z_{\mathcal{U}}
\to \mathcal{U}\subset G(3,\bigwedge\limits^2 V_6)$ is generically finite of degree 
$2$.  The fibers of this morphism over an element of type a,b,c in Lemma
\ref{planecubi} is  respectively in $G(3,V_6\dual)$: $2$ points, $1$ point, and
a rational cubic curve.
\end{proposition}
\proof
Let $(\mu,\nu)$ be an element of $G(3,V_6\dual)\times G(3,\bigwedge\limits^2
V_6)$.  The fiber of a vector bundle at
$\mu$ (resp. $\nu$) will be noted by the name of the bundle with the index $\mu$
(resp. $\nu$). The vector space $K_{3,\mu}$ is isotropic for all the skew-symmetric
forms defined by the elements of $R_{3,\nu}$ if and only if $(\mu,\nu)\in Z$, but
also if and only if the composition:
$$
\begin{CD}
   R_{3,\nu} @>>> \bigwedge\limits^2 V_6  @>>> \bigwedge\limits^2 K_{3,\mu}\dual  
 \end{CD}
$$
is the zero map. So $(\mu,\nu)\in Z \iff R_{3,\nu} \subset E_{12,\mu}$ and we
have the equality $Z=G(3,E_{12})$.

The end of the assertion follows immediatly from  Lemma \ref{planecubi} and
 Remark \ref{3casplans}. \fin
\begin{corollary}
  The locus $\mathcal{U}_c$ in  $\mathcal{U}\subset G(3,\bigwedge\limits^2 V_6)$  of planes of
  type c) has codimension $3$. Consider the following relation  on $\mathcal{U}_c$: $
  p \mathcal{R} p'$ if and only if $p_1(p_2^{-1}(p))=p_1(p_2^{-1}(p'))$. For any
  element $p$ of $ \mathcal{U}_c$, there is a six dimensional subspace $L_p$ of
  $\bigwedge\limits^2 V_6$ such that   the equivalence class of $p$ for
  $\mathcal{R}$ is an open set of $G(3,L_p)$.
\end{corollary}
\proof From the proposition \ref{propisoE12}, for any $p$ in $\mathcal{U}_c$,
$p_1(p_2^{-1}(p))$ is a smooth rational cubic curve $C_p$  in
$G(3,V_6\dual)$. So the restriction of
$E_{12}$ to $C_p$ is $6\OO_{\PP_1}\oplus 6\OO_{\PP_1}(-1)$, and this bundle has
a natural trivial subbundle of rank $6$. Let $L_p$ be the six
dimensional vector subspace of $\bigwedge\limits^2 V_6$ obtained from the image
of this subbundle by the injection of the sequence (\ref{eqE12}).

Proposition \ref{propisoE12} describes $p_1^{-1}(C_p)$ as
the relative Grassmannian $G(3,E_{12|C_p})$. Let $F$ be a
subvector bundle of rank $3$ of $E_{12|C_p}=L_p\otimes\OO_{\PP_1}\oplus
6\OO_{\PP_1}(-1)$. The case c) appears when the line
bundle $\wedge^3 F\dual$ contracts the curve $C_p$. But $\wedge^3 F\dual$ is not
ample if and only if $F$ is a trivial subbundle 
of $L_p\otimes\OO_{\PP_1}$. So $p_1^{-1}(C_p)\cap p_2^{-1}(\mathcal{U}_c)$ is
$(\mathcal{U}\cap G(3,L_p))\times C_p$, and the equivalence classe of $p$ for
$\mathcal{R}$ is $\mathcal{U}\cap G(3,L_p)$. So the
dimension of $\mathcal{U}_c$ is the sum of the dimension of $G(3,6)$ with the
dimension of the family of  rational cubic curves in
$G(3,V_6\dual)$ . In conclusion  $\mathcal{U}_c$ has dimension $33$ and codimension $3$ in
$G(3,\bigwedge\limits^2 V_6)$. \fin
%\begin{proposition}
%  ADD??: Description of the ramification as a determinant of bundles on $Z$
%\end{proposition}
%\subsection{Exemple CY} 

\section{Palatini threefolds}\label{secn4}
In this section we  study the case $n=4$. 
\subsection{D\'efinition and classical properties}
\begin{definition}\label{palatini}
  A smooth $3$ dimensional sub-variety $X$ of $\P5$ is called a Palatini
  threefold\footnote{or a Palatini scroll} if there exists an element of
  $\alpha \in \bigwedge\limits^2 V_6 \otimes W_4$ such that $X=\Proj{F}$ where $F$ is the
  pfaffian vector bundle defined from $\alpha$ in the Definition \ref{defpfaf}
  with $n=4$. It is also classically called (cf \cite{Do} p589) the singular variety
  of the linear system $|W_4\dual|$ of linear line complex in $|V_6\dual|$.
%  by the exact sequence associated to $\alpha$:  
%  \begin{equation}
%    \begin{CD}
%     0@>>>  V_6\dual \otimes \OO_{\pi_{3}}(-1) @>M_\alpha>> V_6 \otimes \OO_{\pi_{3}} @>>> F@>>> 0 
%   \end{CD}
%  \end{equation}
\end{definition}
\begin{notation}\label{notapala}
In this section, denote by $X$  a Palatini threefold in $\P5$, by $h$ the class
of an 
hyperplane of $\P5$, by $S$ the pfaffian cubic surface in $\Pw$ and by $s$
the pullback on $X$ of the class of an hyperplane of $\Pw$. The cotangent bundle
of $\P5$ will be noted $\Omega^1_{\P5}$. 
\end{notation}

So we can immediately obtain the well known resolution of its ideal:
\begin{remark}\label{resolpala}
  The ideal $I_X$ of a Palatini threefold $X$  in $\P5$ has the following
  resolution:
  \begin{equation}
  \label{resopalaP5}
   \begin{CD}
   0@>>>  W_4\dual \otimes \OO_{\P5} @>\alpha>> \Omega^1_{\P5}(2h) @>>> I_X(4h) @>>> 0 
 \end{CD}
\end{equation}
\end{remark}
and the famous equality: $$h^0\OO_X(2h)=h^0\OO_{\P5}(2h)+1.$$

\bigskip
To explain the natural embedding of $X$ in the point/plane incidence  of $\P5$,
F. Zak introduced the following vector bundle:
\begin{definition}
  The canonical extension on $\P5$ displayed in the second column of the
  following diagram of exact sequences
$$
\begin{CD}
  @.  @. 0 @.    @.\\
  @.  @. @VVV @. \\
    0 @>>> W_4\dual \otimes \OO_{\P5}(-h)  @>>> \Omega^1_{\P5}(h) @>>> I_X(3h) @>>>0\\
     @.  @V{\sim}VV @VVV  \\
      0@>>>  W_4\dual \otimes \OO_{\P5}(-h) @>\alpha>> V_6\otimes \OO_{\P5}  \\
      @.  @. @VVV @. \\
      @.  @. \OO_{\P5}(h) @.   @.\\
      @.  @. @VVV @. \\
      @.  @. 0 @.    @.\\
    \end{CD}
$$
  induces on a Palatini threefold $X$ the following extension with middle term a
  rank $3$ vector bundle $E_X$.
  $$
  \begin{CD}
    0 @>>> N_X\dual(3h) @>>> E_X @>>> \OO_X(h) @>>> 0.
  \end{CD}
  $$
  Moreover, the restriction to $X$ of the second line of the previous diagram
  gives the exact sequence:
  \begin{equation}
    \label{resEX}
     \begin{CD}
      0@>>>  \OO_X(-h-s)@>>> W_4\dual \otimes \OO_{X}(-h) @>\alpha>> V_6\otimes \OO_{X} @>>> E_X @>>> 0
     \end{CD}
  \end{equation}
and the determinant of $E_X$ is $\OO_X(3h-s)$.
\end{definition}
\bigskip
From the inclusion $W_4\dual \subset \bigwedge\limits^2 V_6$ and the identification
$W_4=\wedge^3W_4\dual$, we can  consider $\Pwv$  as a subvariety of $G(3,\bigwedge\limits^2 V_6)$. 
\begin{proposition}\label{Z4isoX}
  Let $Z_4$ be the restriction of the isotropic incidence $Z\subset G(3,V_6\dual)\times
  G(3,\bigwedge\limits^2 V_6)$  to  $G(3,V_6\dual)\times \Pwv$. Then $Z_4$ is
  isomorphic to $X$ and the projection from $Z_4$ to $G(3,V_6\dual)$ is the
  natural embedding of $X$ given by $\bigwedge\limits^3 E_X$.
\end{proposition}
\proof Let's first recall the classical description of quadrisecant lines to
$X$. Let $A\dual$ and $B$ be the $3$ dimensional vector subspaces of $V_6\dual$ and
$W_4\dual$ corresponding to a point of $Z_4$. Denote by $A'$ the kernel of 
the surjection from $V_6$ to $A$ and $\PP(A\dual)\subset \P5$ by $\pi_A$. The
restriction of $\Omega^1_{\P5}(1)$ to $\pi_A$ is $A'\otimes \OO_{\pi_A} \oplus
\Omega^1_{\pi_A}(1)$.
% $$B \subset W_4\dual,\ 0\to A' \to V_6 \to A \to 0$$

From the isotropicity of $\pi_A$ with respect to all the
elements of $B$ we see that the restriction of $\alpha$ to $\pi_A$  is the
direct sum of the following maps: 
\begin{center}
$ B \otimes \OO_{\pi_A}(-1) \to A' \otimes\OO_{\pi_A}$ and $ \frac{W_4\dual}{B} \otimes
\OO_{\pi_A}(-1) \to \Omega^1_{\pi_A}(1)$. 
\end{center}
The determinant of the first one gives a cubic curve in $\pi_A\cap X$, and the
second map vanishes on a single (residual) point $\mu$ of $\pi_A\cap X$. So we
have constructed a morphism from $Z_4$ to $X$: $(A\dual,B) \mapsto \mu$.

Moreover, this vanishing shows by  specialization of  sequence
(\ref{resEX}) at the point $\mu$ that the fiber of $E_X\dual$ at $\mu$ is
$A\dual$. So $Z_4$ and $X$ have the same image in $G(3,V_6\dual)$, and the proof
of the proposition is reduced to the proof of the embedding of $Z_4$ to
$G(3,V_6\dual)$. But the fiber of this morphism over the point of
$G(3,V_6\dual)$ corresponding to $A\dual$ is a single point because $A\dual$  is
not isotropic for all the elements of $W_4\dual$. So this projection of $Z_4$ is
one to one, and it must be an isomorphism because the fibers are given by linear
conditions. \fin

\subsection{Anticanonical properties}
Although it is classical that the canonical class $K_X$ of a Palatini threefold
satisfies $K_X^3=-2$ (cf. for instance \cite{Ok}) the following identification and
next proposition seem new.
\begin{lemma}\label{W4W4v}
   The  canonical line
  bundle of $X$:  $\omega_X$ is isomorphic to $\OO_X(s-2h)$. With the notations
  \ref{notapala}, we have from the equality $W_4= H^0(\OO_S(1))$  a canonical
  isomorphism: $$H^0(\omega_X\dual) = W_4\dual $$  
\end{lemma}
\proof The isomorphism $\omega_X\simeq\OO_X(s-2h)$ can be computed directly from
the definition \ref{palatini}. We obtain the isomorphism $H^0(\omega_X\dual) =
W_4\dual $ from the isomorphism between $X$ and $Z_4$ found in the proposition
\ref{Z4isoX} and the fact that $\omega_{Z_4}\dual$ is the pull back of
$\OO_{\Pwv}(1)$. \fin

\begin{proposition}\label{propantican}
  The linear system $|\omega_X\dual|$ has no base points and gives a morphism of
  degree $2$: $$ \begin{CD} X @>2:1>> \Pwv \subset G(3,\bigwedge\limits^2 V_6
    )\end{CD} $$
  %sans lieu base + tres ample
  The anticanonical linear system of $X$ contracts $5$ rational curves. These
  curves have  degree $3$ for the  embeding of
  $X$ in  $\P5$ and also for the embeding of $X$ in $G(3,V_6\dual)$. 
\end{proposition}
\proof The contracted curves of this morphism correspond to the case c) of lemma
\ref{3casplans}: the planes of a Segre. So they are smooth rational cubic curves in
$G(3,V_6\dual)$.  By definition, on such a curve, de divisors $2h$ and  $s$ are
equivalent because $\omega_X\dual=\OO_X(2h-s)$. So those curves have the same
degree with respect to $h$ than to $3h-s$. So the proof will end after the
following:
\def\Fno{\bar{F}}
\begin{lemma}
Let $\Fno$ be the normalized bundle $F(-1)$. The vector space
$$H=H^1\left((S^2\Fno)(-1)\right)$$
has dimension $5$ and it is the kernel
of the following map given by the pfaffians of size $4$  of $M$:
\begin{equation}
  \label{eq:kerpfaf44}
\begin{CD}
0 @>>> H @>>> \bigwedge\limits^2 V_6 = \bigwedge\limits^4 V_6\dual @>>> S^2 W_4 @>>>0  
\end{CD}
\end{equation}
Moreover the ideal of the exceptional locus in
$\Pwv$ of the projection $X=Z_4 \to \Pwv$ is given by 
the $4\times 4$ pfaffians of a skew-symmetric map:
$$\begin{CD}
H \otimes \OO_{\Pwv} (-1) @>>>  H\dual \otimes \OO_{\Pwv} 
\end{CD}.  
$$
\end{lemma}
\proof
Let   $i$ be an isomorphism: $\wedge^2\Fno \to \OO_S$. The restriction of $F$ to a
plane $P$ is of type  c) in  lemma
\ref{3casplans}  if and only if we have $h^1(S^2(\Fno_P))=3$.

To globalize this condition, let's consider the complex:
$$C^\bullet:
\begin{CD}
  0 @>>> V_6\dual\otimes\OO_{\Pw}(-2) @>M>> V_6\otimes\OO_{\Pw}(-1) @>>> 0
\end{CD}.
$$
It is exact in degree $-1$ with cohomology $\Fno$ in degree $0$. The exterior
power of $C^\bullet$ tensorized by $\OO_{\Pw}(2)$ is:
$$\begin{CD}
    0 @>>> S^2 V_6\dual\otimes\OO_{\Pw}(-2) @>>> V_6\dual\otimes
    V_6\otimes\OO_{\Pw}(-1) @>>> \bigwedge\limits^2V_6\otimes\OO_{\Pw} @>>> 0
\end{CD}
$$
with cohomology in degree $(-2,-1,0)$:
$(0,S^2(\Fno)(-1),(\wedge^2(\Fno)(2)))$. So the hyper-cohomology's spectral
sequence of this complex gives the exact sequence (\ref{eq:kerpfaf44}),  the
dimension of $H$, and the vanishings
$$h^0(S^2(\Fno)(-1))=h^2(S^2(\Fno)(-1))=h^0(S^2(\Fno))=h^2(S^2(\Fno))=0.$$

% Now take the symmetric
% power of $C^\bullet$ :
% $$K^{\bullet}:
% \begin{CD}
%     0 @>>> \bigwedge^2 V_6\dual\otimes\OO_{\Pw}(-4) @>>> V_6\dual\otimes
%     V_6\otimes\OO_{\Pw}(-3) @>>> S^2V_6\otimes\OO_{\Pw}(-2) @>>> 0
% \end{CD}.
% $$
% Its cohomology in degree $(-2,-1,0)$ is $(0,(\wedge^2(\Fno))(-3),S^2(\Fno))$.

Now consider the point/plane incidence variety $I_3\subset \Pwv\times \Pw$ and denote by
$p\dual_{3}$ and $p_3$ the first and second projections of this product. We have
the exact sequence:
$$
\begin{CD}
  0@>>> \OO_{\Pwv}(-1)\boxtimes S^2\Fno(-1) @>>> \OO_{\Pwv}\boxtimes S^2\Fno @>>> p_3^*(S^2\Fno) @>>> 0
\end{CD}
$$
From the Leray's spectral sequence and the above vanishings, we have the exact sequence: 
{\small\minCDarrowwidth1em
$$
\begin{CD}
  0@>>> p\dual_{3*}(p_3^*(S^2\Fno)) @>>> H^1(S^2(\Fno)(-1))\otimes
  \OO_{\Pwv}(-1)@>d_M>> H^1(S^2(\Fno))\otimes \OO_{\Pwv} @>>>
  R^1p\dual_{3*}(p_3^*(S^2\Fno)) @>>> 0
\end{CD}
$$
}

Let's now explain how to consider the map $d_M$ as a skew-symmetric map. The isomorphism $i$ gives a symmetric isomorphism $i': S_2(\Fno) \to
S_2(\Fno\dual)$ so the following square is commutative:
$$
\begin{CD}
   (S_2\Fno)(-1) \otimes  S_2\Fno @>i'\otimes id>>  (S_2\Fno\dual)(-1)
   \otimes  S_2\Fno\\
   @Vid \otimes i'VV @VV{\tau}V \\\
  (S_2\Fno)(-1) \otimes  S_2\Fno\dual @>{\tau'}>> \OO_S(-1) 
\end{CD}
$$
The cup-product $H^1\left((S_2\Fno)(-1)\right)\otimes H^1\left((S_2\Fno)(-1)\right) \to
H^2\left((S_2\Fno\otimes S_2\Fno)(-2)\right)$ is anti-commutative, so for any
$z\in W_4$ the
following square is also anti-commutative:
$$
\begin{CD}
   H^1\left((S_2\Fno)(-1)\right) \otimes  H^1\left((S_2\Fno)(-1)\right) @>d_{M,z}\otimes id>>  H^1\left(S_2\Fno\right) \otimes  H^1\left((S_2\Fno)(-1)\right)\\
   @Vid \otimes d_{M,z}VV @VV{\cup}V \\
H^1\left((S_2\Fno)(-1)\right) \otimes  H^1\left(S_2\Fno\right)   @. H^2\left(S_2\Fno\otimes  (S_2\Fno)(-1)\right) \\
   @V{\cup}VV @VV{\overline{\tau\circ(i'\otimes id)}}V \\
H^2\left((S_2\Fno)(-1) \otimes  S_2\Fno\right)   @>\overline{\tau'\circ(id\otimes i')}>> H^2(\OO_S(-1)) 
\end{CD}
$$
In conclusion, the composition:
$$
\begin{CD}
  H \otimes \OO_{\Pwv}(-1)@>d_M>>H^1(S_2\Fno)\otimes \OO_{\Pwv}@>\bar{i'}>>H^1(S_2(\Fno\dual)\otimes \OO_{\Pwv})@>Serre's\ duality>>   H\dual\otimes \OO_{\Pwv}
\end{CD}
$$
is  skew-symmetric and the lemma is proved. Indeed, the  type c) cases 
correspond to the locus where this map has rank at most $2$. \fin

%%%%%%%%%%%%%%%%%%%%%%%%%%%%%%%%

\begin{definition}\label{defPHI1}
 Let $\Sigma_5$ be the symmetric product of order $5$ of $\Pwv$. We define
 the rational map $\Phi_1$ to be:  $$ \Phi_1:
 \begin{array}[t]{ccc}
(W_4\otimes \bigwedge\limits^2
V_6)^{sm}/GL(V_6)  & \dashrightarrow  & \PP(S^3(W_4))\times \Sigma_5\\
\alpha & \mapsto & (S,(h_0\dots h_4))
\end{array}
$$
where $S$ is the pfaffian cubic surface defined by $\alpha$, and 
$h_0,\dots,h_4$ are the five linear sections of $S$  defined in  proposition \ref{propantican}.
\end{definition}
In section \ref{W2V6} we will understand the image of this map.

%%%%%%%%%%%%%%%%%%%%%%%%%%%%%%%%%
%\bigskip
%The following remark 
%\begin{remark}
% T->W quartique marquee et p_*(O_X(h)) pour decrire l'involution 
%\end{remark}

\subsection{Palatini threefolds and endomorphisms}
%(FACULTATIF)
%\begin{remark}
%The incidence $\left\{(p',p)\in X\times X | p\in C_{p'}\right\}$ is the pull back of the
%(point/plane) incidence of $\Pwv \times \Pw$ by the product of the morphisms: $
%|\omega_X\dual|: X \to \Pwv $ and $X \to S \subset \Pw$.
%\end{remark}
Although this part is not required by the main theorem, let's briefly describe
here some connected remarks.

The exceptional geometric properties of a Palatini threefold are classically
considered as natural generalizations of what happens to a Veronese
$\mathcal{V}$ surface embeded in $\PP_4$. Note for instance, in the Veronese situation, the sequence
\ref{resopalaP5} is replaced by:
$$
\begin{CD}
   0@>>>  W_3\dual \otimes \OO_{\PP_4} @>\alpha>> \Omega^1_{\PP_4}(2h) @>>> I_{\mathcal{V}}(3h) @>>> 0 
 \end{CD}.
$$
But the main difference is that in the theory of Severi varieties the embedding
of $\mathcal{V}$ by the complete linear system $|\OO_{\mathcal{V}}(h)|$ is
understood from  an interpretation in terms of matrices of size $3\times 3$ of
rank $1$.    For a Palatini threefold, there is no similar result
to describe the embedding by the complete linear system $|\OO_X(2h)|$. The
following remark  could be a first step in this direction:
\begin{remark}\label{endopala}
The restriction of the line bundle $\omega_X\dual \boxtimes \OO_X(s)$ to the diagonal of $X\times X$ gives the natural
inclusions: $$W_4\dual \otimes W_4 \subset H^0(\OO_X(2h))$$  
In other words, the embedding of a Palatini threefold $X$ with $|\OO_X(2h)|$ has
a canonical projection in $\PP(W_4\dual \otimes W_4)$, and the image of $X$ by
this projection is included in the endomorphisms of $W_4$ of rank $1$.
\end{remark}
\proof It's straighforward from  lemma \ref{W4W4v}. \fin
\section{Geometry in $\bigwedge\limits^2 V_6$}\label{W2V6}
\subsection{Projections from linear spaces}
The Grassmannian variety $G(2,6)$ is one of the $4$ Severi varieties. It is  well
known to have the exceptional property that its projection from a general line
has a unique triple point (cf \cite{I-M}, \cite{Z}). Here, we prove that it has
the same  property with projection from  general $\PP_3$ and points of multiplicity $5$:
\begin{proposition}\label{P4-5sec}
  Denote by
  $\mathcal{U}_5$ the subspace of $G(5,\bigwedge\limits^2 V_6)$ defined by the
  five dimensional vector spaces  such that
  %of $W_4\dual\subset W_5\dual\subset\bigwedge\limits^2 V_6$
   the intersection of their projectivisation  with $G(2,V_6)$ is 
  $5$ linearly independent distinct points. Let $W_4\dual$ be a general $4$-dimensional subspace
  of $\bigwedge\limits^2 V_6$,  then there is a unique element of
  $\mathcal{U}_5$ containing $W_4\dual$.
  % the pfaffian map $\bigwedge\limits^4 V_6\dual
  % \to S_2 W_5$ has rank $10$. 

  % Let $N_{W_5}$ be
  % the vector space spanned by elements of $W_5\dual$ of rank at most $2$. Then
  % $N_{W_5}$  is the kernel of the pfaffian map:
  % $$
  % \begin{CD}
  %   0 @>>> N @>>> \bigwedge\limits^2 V_6\simeq \bigwedge\limits^4 V_6\dual @>Pf>> S_2 W_5 @>>>0 
  % \end{CD}
  % $$
\end{proposition}
\proof
First remark that the incidence variety
$$I_{4,5}=\{(W_4\dual,W_5\dual)| W_4\dual\subset W_5\dual
\subset \bigwedge\limits^2 V_6,
W_5\dual \in \mathcal{U}_5\}$$ has the same dimension as
$G(4,\bigwedge\limits^2V_6)$, so we have to prove that the natural projection is
birational. 

So,  consider a general element $W_4\dual$ in the image of this projection, and
chose  an element $W_5\dual$ such that $(W_4\dual,W_5\dual)\in I_{4,5}$. Denote by $\Pw,\pi_4$ their projectivisation. 
The vector space $H^0(I_{\Pw \cup G(2,V_6)}(2))$ is the kernel of the 
map $\bigwedge\limits^4 V_6\dual\to S_2 W_4$. So it has dimension $5$. Now
remark that we also have $h^0(I_{\pi_4 \cup  G(2,V_6)}(2))=5$ because the ideal
of the $5$ points $\pi_4 \cap G(2,V_6)$ in $\pi_4$ is
a $10$ dimensional space of quadrics. So we proved that $\pi_4$ must be in all
the quadrics of $H^0(I_{\Pw \cup G(2,V_6)}(2))$. It gives the following
linear conditions satisfied by any $W_5\dual$ of $\mathcal{U}_5$ containing $W_4\dual$:
$$ W_5\dual \subset \bigcap\limits_{q \in H^0(I_{\Pw \cup G(2,V_6)}(2)) } (W_4\dual)^{\bot_q}$$
where $\bot_q$ denotes the orthogonal with respect to the
quadratic form $q$ on $\bigwedge^2 V_6$. So  unicity of $W_5\dual$
will be a corollary of  existence of an exemple of $W_4$ such that
$\bigcap\limits_{q \in H^0(I_{\Pw \cup G(2,V_6)}(2)) } (W_4\dual)^{\bot_q}$ has
dimension $5$ as it is the case in the following:
\begin{example}
Let's consider a basis $(\eps_i)$ of $V_6$, and the $5$ elements  $$u_0=\eps_0\wedge
\eps_3, u_1=\eps_1 \wedge \eps_4, u_2=\eps_2 \wedge \eps_5, u_3=(\eps_0+\eps_1+\eps_2)\wedge (\eps_4+\eps_3+\eps_5),u_4=(\eps_1+\eps_4+\eps_2)\wedge (\eps_3+\eps_1+\eps_5).$$
Denote by $W_5\dual$ the $5$ dimensional vector space spanned by the $(u_i)$ and
$$W_4\dual=\{\sum_{0\leq i \leq 4} \lambda_i.u_i | \sum_{0\leq i \leq 4} \lambda_i =0\}.$$
Then $\bigcap\limits_{q \in H^0(I_{\Pw \cup G(2,V_6)}(2)) } (W_4\dual)^{\bot_q}$
has dimension $5$.
\end{example}
\proof
%lll1=(e_0+e_1+e_2)*(e_4+e_3+e_5)
%lll2=(e_1+e_4+e_2)*(e_3+e_1+e_5)
%g4=sum flatten for i from 0 to 5 list for j from i+1 to 5 list coefficient(e_i*e_j,lll1)*p_(i,j)
%g5=sum flatten for i from 0 to 5 list for j from i+1 to 5 list coefficient(e_i*e_j,lll2)*p_(i,j)
%tmpP4=gens ideal(p_(0,3),p_(1,4),p_(2,5),g4,g5)
We can compute with \cite{Macaulay2} that  $H^0(I_{\Pw\cup G(2,V_6)}(2))$ is generated
by the following five quadrics in Plucker coordinates:
\begin{itemize}
\item[.] ${p}_{({3},{4})} {p}_{(1,{5})}-{p}_{(1,{4})} {p}_{({3},{5})}+{p}_{(1,{3})}
  {p}_{({4},{5})}$
  
\item[.]$  {p}_{(1,{2})} {p}_{(0,{5})}-{p}_{({2},{4})} {p}_{(0,{5})}-{p}_{(0,{2})}
{p}_{(1,{5})}+{p}_{({2},{3})} {p}_{(1,{5})}+{p}_{(0,1)} {p}_{({2},{5})}-{p}_{(1,{3})}
{p}_{({2},{5})}+{p}_{(0,{4})} {p}_{({2},{5})}-{p}_{({3},{4})}
{p}_{({2},{5})}+{p}_{(1,{2})} {p}_{({3},{5})}+{p}_{({2},{4})}
{p}_{({3},{5})}-{p}_{(0,{2})} {p}_{({4},{5})}-{p}_{({2},{3})} {p}_{({4},{5})}$

\item[.]  $     {p}_{({2},{3})} {p}_{(0,{4})}-{p}_{(0,{3})} {p}_{({2},{4})}+{p}_{(0,{2})}
{p}_{({3},{4})}-{p}_{(1,{3})} {p}_{(0,{5})}+{p}_{({2},{4})} {p}_{(0,{5})}-{p}_{({3},{4})}
{p}_{(0,{5})}+{p}_{(0,{3})} {p}_{(1,{5})}-{p}_{({2},{3})} {p}_{(1,{5})}+{p}_{(1,{3})}
{p}_{({2},{5})}-{p}_{(0,{4})} {p}_{({2},{5})}+{p}_{({3},{4})} {p}_{({2},{5})}-{p}_{(0,1)}
{p}_{({3},{5})}-{p}_{(1,{2})} {p}_{({3},{5})}+{p}_{(0,{4})}
{p}_{({3},{5})}-{p}_{({2},{4})} {p}_{({3},{5})} +{p}_{(0,{2})}
{p}_{({4},{5})}-{p}_{(0,{3})} {p}_{({4},{5})}+{p}_{({2},{3})} {p}_{({4},{5})}$

\item[.] $       {p}_{(1,{2})} {p}_{(0,{4})}-{p}_{(0,{2})} {p}_{(1,{4})}+{p}_{(0,1)}
{p}_{({2},{4})}-{p}_{({2},{4})} {p}_{(0,{5})}+{p}_{({2},{3})} {p}_{(1,{5})}-{p}_{(1,{3})}
{p}_{({2},{5})}+{p}_{(0,{4})} {p}_{({2},{5})}-{p}_{({3},{4})}
{p}_{({2},{5})}+{p}_{(1,{2})} {p}_{({3},{5})}+{p}_{({2},{4})}
{p}_{({3},{5})}-{p}_{(0,{2})} {p}_{({4},{5})}-{p}_{({2},{3})} {p}_{({4},{5})}$
\item[.] $       {p}_{(1,{2})} {p}_{(0,{3})}-{p}_{(0,{2})} {p}_{(1,{3})}+{p}_{(0,1)} {p}_{({2},{3})}-{p}_{({2},{4})} {p}_{(0,{5})}+{p}_{({3},{4})} {p}_{(0,{5})}+{p}_{({2},{3})} {p}_{(1,{5})}-{p}_{(1,{3})} {p}_{({2},{5})}+{p}_{(0,{4})} {p}_{({2},{5})}-{p}_{({3},{4})} {p}_{({2},{5})}+{p}_{(1,{2})} {p}_{({3},{5})}-{p}_{(0,{4})} {p}_{({3},{5})}+{p}_{({2},{4})} {p}_{({3},{5})}-{p}_{(0,{2})} {p}_{({4},{5})}+{p}_{(0,{3})} {p}_{({4},{5})}-{p}_{({2},{3})} {p}_{({4},{5})}$
\end{itemize}
%
% toString I5
% --matrix {{p_(3,4)*p_(1,5)-p_(1,4)*p_(3,5)+p_(1,3)*p_(4,5),p_(1,2)*p_(0,5)-p_(2,4)*p_(0,5)-p_(0,2)*p_(1,5)+p_(2,3)*p_(1,5)+p_(0,1)* p_(2,5)-p_(1,3)*p_(2,5)+p_(0,4)*p_(2,5)-p_(3,4)*p_(2,5)+p_(1,2)*p_(3,5)+ p_(2,4)*p_(3,5)-p_(0,2)*p_(4,5)-p_(2,3)*p_(4,5), p_(2,3)*p_(0,4)-p_(0,3)*p_(2,4)+p_(0,2)*p_(3,4)-p_(1,3)*p_(0,5)+p_(2,4)* p_(0,5)-p_(3,4)*p_(0,5)+p_(0,3)*p_(1,5)-p_(2,3)*p_(1,5)+p_(1,3)*p_(2,5)- p_(0,4)*p_(2,5)+p_(3,4)*p_(2,5)-p_(0,1)*p_(3,5)-p_(1,2)*p_(3,5)+p_(0,4)*p_(3,5)-p_(2,4)*p_(3,5)+p_(0,2)*p_(4,5)-p_(0,3)*p_(4,5)+p_(2,3)*p_(4,5),p_(1,2)*p_(0,4)-p_(0,2)*p_(1,4)+p_(0,1)*p_(2,4)-p_(2,4)*p_(0,5)+p_(2,3)*p_(1,5)-p_(1,3)*p_(2,5)+p_(0,4)*p_(2,5)-p_(3,4)*p_(2,5)+p_(1,2)*p_(3,5)+p_(2,4)*p_(3,5)-p_(0,2)*p_(4,5)-p_(2,3)*p_(4,5),p_(1,2)*p_(0,3)-p_(0,2)*p_(1,3)+p_(0,1)*p_(2,3)-p_(2,4)*p_(0,5)+p_(3,4)* p_(0,5)+p_(2,3)*p_(1,5)-p_(1,3)*p_(2,5)+p_(0,4)*p_(2,5)-p_(3,4)*p_(2,5)+ p_(1,2)*p_(3,5)-p_(0,4)*p_(3,5)+p_(2,4)*p_(3,5)-p_(0,2)*p_(4,5)+p_(0,3)* p_(4,5)-p_(2,3)*p_(4,5)}}
%--toString IP4
%--ideal(p_(3,5),p_(0,5)-p_(1,5)+p_(4,5),p_(3,4)+p_(4,5),p_(2,4)-p_(1,5)+p_(4,5),p_(0,4)-p_(1,5)+p_(4,5),p_(2,3)-p_(1,5),p_(1,3)-p_(1,5),p_(1,2)+p_(4,5),p_(0,2),p_(0,1))
and check that the ideal of the orthogonal of $\Pw$ with respect to these $5$
quadrics is generated by the $10$ independant equations:
$       ({p}_{({3},{5})},{p}_{(0,{5})}-{p}_{(1,{5})}+{p}_{({4},{5})},{p}_{({3},{
    4})}+{p}_{({4},{5})},{p}_{({2},{4})}-{p}_{(1,{5})}+{p}_{({4},{5})},{p
}_{(0,{4})}-{p}_{(1,{5})}+{p}_{({4},{5})},{p}_{({2},{3})}-{p}_{(1,{5
  })},{p}_{(1,{3})}-{p}_{(1,{5})},{p}_{(1,{2})}+{p}_{({4},{5})},{p}_{(0,{2
  })},{p}_{(0,1)})$. So this example completes the  proof of the birationality of the projection
from $I_{4,5}$ to $G(4,\bigwedge\limits^2 V_6)$. So we have proved  proposition
\ref{P4-5sec}. \fin

\begin{corollary}\label{defPHI2}
With notations of  definition \ref{inscribed}, we can define 
 the rational map $\Phi_2$ by:  $$ \Phi_2:
 \begin{array}[t]{ccc}
(W_4\otimes \bigwedge\limits^2
V_6)^{sm}/GL(V_6)  & \dashrightarrow  & \mathcal{H}\\
\alpha & \mapsto & (S,(H_0\dots H_4))
\end{array}
$$
where $S$ is the pfaffian cubic surface defined by $\alpha$, 
and $H_i$ is defined like this:

From proposition \ref{P4-5sec}, consider  the five
points $(u_i)_{0\leq i \leq 4}$  of $G(2,V_6)$ such that $\Pw$ is in the linear
span $<(u_i)_{0\leq i \leq 4}>$. Then take:
$$H_i=\Pw\cap <(u_j)_{0\leq j\leq 4, j\neq  i} >$$
\end{corollary}
\proof After proposition \ref{P4-5sec}, we only have to explain why
$H_0,\dots,H_4$ is inscribed on $S$. But for $\{i_0,\dots,i_4\}=\{0,\dots,4\}$
the point $H_{i_0}\cap H_{i_1}\cap H_{i_2}$ is on the line $(u_{i_3},u_{i_4})$
so it corresponds to a matrix of rank $4$ and  is on $S$. \fin

\begin{remark}\label{Hrational}
  The variety $\mathcal{H}$ is rational of dimension $24$.
\end{remark}
\proof Let $\Sigma_5'$ be the image of $\mathcal{H}$ in $|\OO_{\Pw}(5)|$ by the
second projection.  It is an openset of the symmetric product $\Sigma_5$
defined in \ref{defPHI1}. So it is a rational $15$-dimensional variety (cf \cite{GKZ} Th 2.8 p 137).

The partial derivatives of order $2$ of any element of
$\Sigma_5'$ are linearly independent cubic forms. So they give a rank $10$ subsheaf $\mathcal{F}_2$ of
$H^0(\OO_{\Pw}(3))\otimes \OO_{\Sigma_5'}$ locally free with respect to Zariski's topology.

Now remark that $\mathcal{H}$ is the openset  of $\PP(\mathcal{F}_2)$
corresponding to smooth cubic surfaces. So $\mathcal{H}$ is rational of
dimension $24$. \fin

\subsection{An explicit formula and  proof of  theorem \ref{mainbir}}\label{sectionexplicit}

Surprisingly, we are able to give in this section an explicit formula.
Recently, a explicit result was also found by
F. Tanturri in   \cite{T}: An algorithm to obtain a pfaffian
representation from a cubic equation.  The two main difference, are the
following:

-first he wants to find any pfaffian representation of $S$, but here we need to
find a unique point in the moduli space.

-The construction starts with five points on $S$, so it is a problem of extending
the $5$ by $5$ skew-symmetric matrix of the resolution of the $5$ points to a
$6$ by $6$ one with pfaffian $S$, while we start with an inscribed pentahedron.

% \begin{remark}\label{equation}
%   The element $(S,\Pi)$ is in $\mathcal{H}$ if and only if the equation of $S$
%   can be written:
%   $$ \sum_{0\leq i<j<k\leq 4} A_{i,j,k}.h_i.h_j.h_k =0 ,  h_4=\sum_{0\leq i\leq 3} h_i $$
%   where $(h_i)_{0\leq i \leq 4}$ are  equations of the elements of $\Pi$.
% \end{remark}
\begin{lemma}\label{equationnorm}
  Let $(x_i)_{0\leq i \leq 3}$ be a basis of $W_4$, and $\mathcal{A}_9$ be the
  following subspace of $\complex^{10}\times \PP_4$:
  $$
 \mathcal{A}_9=\left\{\left((a_{i,j,k})_{0\leq i<j<k\leq 4},(b_i)_{0\leq i\leq 4} \right) \left|
   \begin{array}[c]{l}
      a_{0,1,4}=1\mbox{ and for } 0\leq i \leq 4,\ b_i\neq 0, \\
     \mbox{and for }  0\leq i<j<k\leq 3,\ a_{i,j,k}=1     
   \end{array}\right.
 \right\}.$$
 Then the following map is birational:
 \begin{equation}
   \begin{array}[c]{ccc}
     PGL_4\times \mathcal{A}_9 & \to & \mathcal{H}_{ord}\\
     (P,\left((a_{i,j,k})_{0\leq i<j<k\leq 4},(b_i)_{0\leq i\leq 4} \right)) &
     \mapsto & (S,(H_0,\dots,H_4))
   \end{array}
 \end{equation}
 where $$\sum_{0\leq i<j<k\leq 4} a_{i,j,k}.w_i.w_j.w_k =0,
 $$ is an equation of $S$, and for all $0\leq i \leq 4$, $w_i=0$ is an equation of
 $H_i$ with the following equalities:
 $w_4=\sum\limits_{i=0}^3 \frac{b_4.w_i}{b_i}$, {\footnotesize$\left(\begin{array}[c]{l}
   w_0\\
   w_1\\
   w_2\\
   w_3
 \end{array}\right)=P. \left(\begin{array}[c]{l}
   x_0\\
   x_1\\
   x_2\\
   x_3
 \end{array}\right)
$}.
\end{lemma}
\proof
Let $\Pi=(H_0,\dots,H_4)$ be an ordered pentahedron and $P'$ be the unique projective transformation that sends the ordered pentahedron
$(x_0,x_1,x_2,x_3,x_0+x_1+x_2+x_3)$ to $(H_0,\dots,H_4)$. Denote by $h_i$ the
equation of $H_i$ defined by: {\footnotesize$\left(\begin{array}[c]{l}
   h_0\\
   h_1\\
   h_2\\
   h_3
 \end{array}\right)=P'. \left(\begin{array}[c]{l}
   x_0\\
   x_1\\
   x_2\\
   x_3
 \end{array}\right)
$} and $h_4=\sum_{i=0}^3 h_i$. Cubic surfaces $S$ such that $(S,\Pi)$ is in
$\mathcal{H}_{ord}$ are the smooth surfaces with equation:
  $$ \sum_{0\leq i<j<k\leq 4} A_{i,j,k}.h_i.h_j.h_k =0,\  (A_{i,j,k})_{(0\leq i < j
  < k \leq 4)} \in \PP_9.$$
Now remark that the map:
$$
\begin{array}[c]{lll}
  \mathcal{A}_9 & \to \PP_9\\
  (a,b) & \mapsto & (A_{i,j,k}=a_{i,j,k}.b_i.b_j.b_k)_{0\leq i <j <k \leq 4}
\end{array}
$$
is birational because we can compute its inverse with the following formulas\footnote{If
  one works with affine spaces instead of $\PP_9$ and $\PP_4$, then one  needs to
  extract a cubic root to solve the equalities.}:
$$
\frac{b_0}{b_3}=\frac{A_{0,1,2}}{A_{1,2,3}},\frac{b_1}{b_3}=\frac{A_{0,1,2}}{A_{0,2,3}},\frac{b_2}{b_3}=\frac{A_{0,1,2}}{A_{0,1,3}},\frac{b_4}{b_3}=\frac{A_{0,1,4}}{A_{0,1,3}},
a_{i,j,4}=\frac{A_{i,j,4}}{A_{i,j,3}}\cdot \frac{A_{0,1,3}}{A_{0,1,4}}.$$
So we obtain the lemma from the equalities: $0\leq i \leq 4, w_i=h_i.b_i$ with
$P$ defined by the product of the diagonal matrix
$(\frac{b_0}{b_4},\dots,\frac{b_3}{b_4})$ with $P'$. \fin

\begin{definition}\label{defexplicit}
  Let $\mathcal{A}'_9$ be the set of triples $(a,b,u)$  such that  $(a,b)$ is an
  element of $\mathcal{A}_9$ defining a smooth cubic surface:
$$\sum_{0\leq i<j<k\leq 4} a_{i,j,k}.w_i.w_j.w_k =0,\ w_4=\sum\limits_{i=0}^3 \frac{b_4.w_i}{b_i},$$
and $u$ is a root of the following equation in $X$: $$X^{2}+X\cdot
  (1+a_{0,2,4}-a_{0,3,4})+a_{0,2,4}=0.$$
Denote by $v=-(1+a_{0,2,4}-a_{0,3,4})-u$ the other one and define:
  $$e_1=a_{0,2,4}+a_{1,2,4}-a_{2,3,4},\ e_2=1+a_{1,2,4}-a_{1,3,4},\ e_3=(-a_{1,2,4}+a_{1,3,4}-1) v-a_{1,2,4}-a_{0,2,4}+a_{2,3,4} $$
  $$M_4=\left(\begin{array}{cccccc}
0 & u & -1 & a_{1,2,4} & e_1 & e_2 \\
-u & 0 & 0 & 0 & a_{0,2,4} & -u \\
1 & 0 & 0 & 0 & -v & 1 \\
-a_{1,2,4} & 0 & 0 & 0 & a_{1,2,4} v & -a_{1,2,4} \\
-e_1 & -a_{0,2,4} & v & -a_{1,2,4} v & 0 & e_3 \\
-e_2 & u & -1 & a_{1,2,4} & -e_3 & 0
\end{array}\right) $$
$$M_{0123}=\left(\begin{array}{cccccc}
0 & 0 & 0 & w_{0}+w_{3} & w_{3} & w_{3} \\
0 & 0 & 0 & w_{3} & w_{1}+w_{3} & w_{3} \\
0 & 0 & 0 & w_{3} & w_{3} & w_{2}+w_{3} \\
-w_{0}-w_{3} & -w_{3} & -w_{3} & 0 & 0 & 0 \\
-w_{3} & -w_{1}-w_{3} & -w_{3} & 0 & 0 & 0 \\
-w_{3} & -w_{3} & -w_{2}-w_{3} & 0 & 0 & 0
\end{array}\right) $$
\end{definition}
% e1:=x02+y12-y23;e2:=1+y12-y13;e3:=(-v-v*y12+v*y13-x02-y12+y23);
% M4:=[[0,u,-1,y12,e1,e2],[0,0,0,0,x02,-u],[0,0,0,0,-v,1],[0,0,0,0,(y12)*v,-y12],[0,0,0,0,0,e3],[0$6]]:;M4:=M4-transpose(M4)
\begin{theorem}\label{explicit}
For a generic element $(P,(a,b,u))$ of $PGL_4\times \mathcal{A}'_9$, the element
$M$ of $(W_4\otimes \bigwedge\limits ^2 V_6)$  defined by
$M=M_{0123}+w_4 M_4$ is such that: $\Phi_1(M)=\Phi_2(M)=(S,\Pi)$ where
the equation of $S$ and $\Pi$ are given by the formula in lemma \ref{equationnorm}.
\end{theorem}
\proof
The difficulty was to find $M_4$.  It was done by tracking the rational cubic curve in
$\P5$ associated to the plane $w_4=0$ in proposition \ref{propantican}.  But
now that we have found $M_4$, it is much easier to check that $M$ safisties the required properties.

NB: To obtain a more compact presentation, we have glued the indexes of the
$a_{i,j,k}$ in the next formulas.

- First, one can check that the pfaffian of $M$ is $$
a_{024} w_{0} w_{2} w_{4}+ a_{034} w_{0} w_{3} w_{4}+a_{234} w_{2}
w_{3} w_{4} +a_{124} w_{1} w_{2} w_{4} +a_{134} w_{1} w_{3} w_{4}  +
w_{0} w_{1} w_{4}+\sum_{0\leq i<j<k\leq 3} w_i w_j w_k$$
% a024*w_0*w_2*w_4+
%w_0*w_2*w_3+
%w_0*w_2*w_1+
%w_0*w_4*w_3*a034+
%w_0*w_4*w_1+
%w_0*w_3*w_1+
%w_2*w_4*w_3*a234+
%w_2*w_4*w_1*a124+
%w_2*w_3*w_1+
%w_4*w_3*w_1*a134

- Now to prove that $\Phi_2(\alpha)=(S,\Pi)$ we just have to remark that $M_4$ has
rank $2$, and also the $4$ values of $M_{0123}$ at the points where
$(w_0,w_1,w_2,w_3)$ take the values $(1,0,0,0)$, $(0,1,0,0)$, $(0,0,1,0)$, $(0,0,0,1)$.

- To obtain that $\Phi_1(\alpha)=(S,\Pi)$ we need to find $5$ elements $(P_i)$
of $GL(V_6)$ such that $\t{P_i}.M.P_i =
\left(\begin{array}[c]{cc}
  0 & A_i\\
  -A_i & 0
\end{array}\right)
$ where $A_i$ are  $3$ by $3$ symmetric matrices with linear entries. We found
the following ones easily,

$$P_4=Id, P_3=\left(\begin{array}{cccccc}
0 & 0 & 0 & 0 & 0 & 1 \\
\frac{v}{u} & 0 & 0 & \frac{(-a_{024}-a_{124}+a_{234})}{u} & 0 & 0 \\
0 & 1 & 0 & 0 & 1+a_{124}-a_{134} & 0 \\
0 & 0 & \frac{-1}{a_{124}} & 0 & 0 & 0 \\
0 & 0 & 0 & 1 & 0 & 0 \\
0 & 0 & 0 & 0 & 1 & 0
\end{array}\right) $$
but the next ones only after understanding that we should use the $SL_2\times SL_2
\times SL_2$ action that preserves the $3$ marked lines in the
intersection of the two Segre $\PP_1\times \PP_2$ defined by $w_i=0$ and $w_4=0$

$$P_1=\left(\begin{array}{cccccc}
0 & 0 & \frac{-a_{024}}{a_{234}} & 0 & 0 & 0 \\
\frac{(-u) (a_{024}+u)}{a_{024} (u+1)} & 1 & \frac{a_{024}}{a_{234}} & \frac{a_{024} a_{134}-a_{024} u-a_{024}+u a_{234}}{a_{024} (u+1)} & \frac{-a_{124}}{u} & 0 \\
\frac{u (a_{024}+u)}{a_{024} (u+1)} & 0 & 0 & \frac{-a_{024} a_{134}+a_{024} u+a_{024}-u a_{234}}{a_{024} (u+1)} & 0 & 0 \\
0 & 0 & \frac{u}{a_{234}} & 0 & 0 & -1 \\
0 & 0 & \frac{-u}{a_{234}} & -1 & 1 & 1 \\
0 & 0 & 0 & 1 & 0 & 0
\end{array}\right) $$

$$P_2=
\left(\begin{array}{cccccc}
\frac{1}{a_{134}} & 0 & 0 & 0 & 0 & 0 \\
0 & \frac{1+v}{u+1} & 0 & 0 & \frac{-v-a_{024}+a_{234}+v a_{134}}{u+1} & 0 \\
-\frac{1}{a_{134}} & \frac{-1-v}{u+1} & 1 & 0 & \frac{v+a_{024}-a_{234}-v a_{134}}{u+1} & a_{124} \\
\frac{1}{a_{134}} & 0 & 0 & 1 & 0 & 0 \\
0 & 0 & 0 & 0 & 1 & 0 \\
-\frac{1}{a_{134}} & 0 & 0 & -1 & -1 & 1
\end{array}\right) 
$$.

$$P_0=
\left(\begin{array}{cccccc}
a_{124} & \frac{-\left(a_{024} u+a_{024}-u a_{234}\right)^{2}}{u^{2} a_{234}} & -a_{134} & \frac{-a_{124} a_{024}}{a_{024} u+a_{024}-u a_{234}} & 0 & \frac{a_{024} a_{134}-a_{024} u-a_{024}+u a_{234}}{a_{024} u+a_{024}-u a_{234}} \\
0 & \frac{\left(a_{024} u+a_{024}-u a_{234}\right)^{2}}{u^{2} a_{234}} & 0 & 0 & 0 & 0 \\
0 & 0 & a_{134} & 0 & 0 & \frac{-a_{024} a_{134}+a_{024} u+a_{024}-u a_{234}}{a_{024} u+a_{024}-u a_{234}} \\
1 & \frac{(u+1) (a_{024} u+a_{024}-u a_{234})}{a_{234} u} & 0 & \frac{u (a_{024}-a_{234})}{a_{024} u+a_{024}-u a_{234}} & -1 & -1 \\
0 & \frac{(-u-1) (a_{024} u+a_{024}-u a_{234})}{a_{234} u} & 0 & 0 & 1 & 0 \\
0 & 0 & 0 & 0 & 0 & 1
\end{array}\right) 
$$
and we have proved  theorem \ref{explicit}. \fin

\bigskip
We are now able to obtain a more explicit version of  theorem \ref{mainbir} stated in  introduction.
\begin{corollary}\label{bir}
  Maps $\Phi_1$ and $\Phi_2$ coincide on a open set, and  give birational maps:
  $$(W_4\otimes \bigwedge\limits^2 V_6)^{sm}/GL(V_6) \dashrightarrow \mathcal{H}.$$
\end{corollary}
\proof 

First remark that both spaces are irreducible of dimension $24$. Now
consider with  notations of  definition \ref{defexplicit} the following map:
\begin{center}
$
\begin{array}[c]{lll}
 PGL_4\times \mathcal{A}'_9 & \to & W_4\otimes \bigwedge\limits^2 V_6\\
 (P,a,b,u) & \mapsto & M_{0123}+w_4.M_4
\end{array}
$, where $w_4=\sum\limits_{i=0}^3
\frac{b_4.w_i}{b_i}$, and {\footnotesize$\left(\begin{array}[c]{l}
   w_0\\
   w_1\\
   w_2\\
   w_3
 \end{array}\right)=P. \left(\begin{array}[c]{l}
   x_0\\
   x_1\\
   x_2\\
   x_3
 \end{array}\right)
$}.
\end{center}
and denote by $f$ its composition with the canonical projection from $
(W_4\otimes \bigwedge\limits^2 V_6)^{sm}$ to $ (W_4\otimes
\bigwedge\limits^2 V_6)^{sm}/GL(V_6)$.

The  map $PGL_4\times\mathcal{A}_9' \to PGL_4\times \mathcal{A}_9$ has degree $2$
because of the permutation of $u$ and $v$, and the rational map $PGL_4\times
\mathcal{A}_9$ to $\mathcal{H}$ has degree $5!$ from the choice of the order and lemma \ref{equationnorm}. So we have from  theorem
\ref{explicit} the commutative diagram of rational maps:
$$
\begin{CD}
  PGL_4\times\mathcal{A}_9' @>2:1>>  PGL_4\times \mathcal{A}_9 @>1:1>>  \mathcal{H}_{ord} \\
   @VVfV  @. @VV(5!):1V\\
  (W_4\otimes \bigwedge\limits^2 V_6)^{sm}/GL(V_6) @.@>\Phi_1>\Phi_2> \mathcal{H} 
\end{CD}
$$
So $\Phi_1$ and $\Phi_2$ are dominant and coincide on an open set, and we
 just have to prove that $f$  has degree $2.(5!)$ also. We will do
this by providing an example of $(S,\Pi)\in \mathcal{H}$ such that the permutation of $u$ and $v$, and the
permutations of the elements of $\Pi$ can be obtained by the action of
$GL(V_6)$. It is more convenient to take an example where all the elements in
the preimage of $(S,\Pi)$  in $PGL_4\times \mathcal{A}'_9$ have all the same values for $(a)$
and $(b)$. So we end the proof with the following invariant example:
\begin{example}(Klein-Sylvester)
  With the following
  values: $u=e^{\frac{2i\pi}{3}}, v=e^{\frac{-2i\pi}{3}}.$ 
  for $0\leq i < j < k \leq 4,\ a_{i,j,k}=1$. The permutation of
  $u$ with $v$, and also the permutations of the $(w_i)_{0\leq i\leq 4}$ can be obtained from
  the action of $GL(V_6)$ on $M=M_{0123}+w_4.M_4$. Note that  if we add the conditions $b_i=-b_4$ for $0
  \leq i \leq 3$, this is the case of the Klein cubic with its Sylvester Pentahedron).
\end{example}
\proof
Denote by $P_T=I_3\otimes
\left(\begin{array}[c]{cc}
  t_0 & t_1\\
  t_2 & t_3
\end{array}\right)
$ the matrix {\footnotesize$\left(\begin{array}{cccccc}
t_{0} & 0 & 0 & t_{1} & 0 & 0 \\
0 & t_{0} & 0 & 0 & t_{1} & 0 \\
0 & 0 & t_{0} & 0 & 0 & t_{1} \\
t_{2} & 0 & 0 & t_{3} & 0 & 0 \\
0 & t_{2} & 0 & 0 & t_{3} & 0 \\
0 & 0 & t_{2} & 0 & 0 & t_{3}
\end{array}\right)  $} and remark that $\t{P_T}M_{0123} P_T = M_{0123}$ when $\left|\begin{array}[c]{cc}
  t_0 & t_1\\
  t_2 & t_3
\end{array}\right|=1$. For a square matrix $\mathcal{T}$, let $D_{\mathcal{T}}$
be the block diagonal matrix
$\left(\begin{array}[c]{cc}
  \mathcal{T}& 0\\
  0 & \mathcal{T}
\end{array}\right)$. So we will first use matrices like $D_\mathcal{T}$  to
obtain the desired form in the plane $w_4=0$ and then correct the last matrix
with $P_T$. We found the following matrices:
\begin{itemize}
\item[.] Permutation of $u$ and $v$:  $P_{uv}=I_3\otimes\left(\begin{array}{cc}
\frac{i.u}{\sqrt{2}} & \frac{\sqrt{6}}{2} \\
\frac{-\sqrt{6}}{2} & \frac{i.v}{\sqrt{2}}
\end{array}\right) $ then $\t{P_{uv}} M_4 P_{uv} = \overline{M_4}.$
\item[.] $ T_{01}=\left(\begin{array}{ccc}
0 & 1 & 0 \\
1 & 0 & 0 \\
0 & 0 & 1
\end{array}\right) $, and
% $P_{01}=I_3\otimes
% \frac{\left(\begin{array}{cc}
% u & v+2 \\
% -1 & -u
% \end{array}\right)}{\sqrt{2}} $
$P_{01}=I_3\otimes
\left(\begin{array}{cc}
\frac{u}{\sqrt{2}} & \frac{v+2}{\sqrt{2}} \\
\frac{-1}{\sqrt{2}} & \frac{-u}{\sqrt{2}}
\end{array}\right)$,  the conjugation $\t{(D_{T_{01}}.P_{01})}.M.(D_{T_{01}}.P_{01})$
    permutes $w_0$ and $w_1$.
\item[.]   $T_{02}=\left(\begin{array}{ccc}
0 & 0 & 1 \\
0 & 1 & 0 \\
1 & 0 & 0
\end{array}\right) $, $P_{02}=I_3 \otimes \left(\begin{array}{cc}
 \frac{e^{\frac{i\pi}{3}}}{\sqrt{2}}& \frac{-i \sqrt{6}}{2} \\
 \frac{e^{\frac{-i\pi}{3}}}{\sqrt{2}}& \frac{v}{\sqrt{2}}
\end{array}\right) $, the conjugation $\t{(D_{T_{02}}.P_{02})}.M.(D_{T_{02}}.P_{02})$ permutes $w_0$ and $w_2$. 
\item[.]   $T_{03}=\left(\begin{array}{ccc}
1 & 1 & 1 \\
0 & -1 & 0 \\
0 & 0 & -1
\end{array}\right) $, 
  $P_{03}=I_3 \otimes \left(\begin{array}{cc}
 \frac{i \sqrt{6}}{2} & \frac{-u}{\sqrt{2}}\\
  \frac{-v}{\sqrt{2}} & \frac{-i \sqrt{6}}{2}
\end{array}\right) $, the conjugation $\t{(D_{T_{03}}.P_{03})}.M.(D_{T_{03}}.P_{03})$ permutes $w_0$ and $w_3$.
\item[.]   $T_{34}=
  \left(\begin{array}{ccc}
0 & -\frac{1}{v} & 0 \\
0 & 0 & 1 \\
1 & 0 & 0
\end{array}\right)$ , $P_{34}=I_3\otimes  \left(\begin{array}{cc}
 -\frac{v}{\sqrt{2}}&  -\frac{i \sqrt{6}}{2}\\
  -\frac{i \sqrt{6}}{2} &\frac{u}{\sqrt{2}}
\end{array}\right)$, then $\t{(P_3 D_{T_{34}} P_{34})}.M.(P_3 D_{T_{34}} P_{34})$ permutes
  $w_4$ and $w_3$ with the matrix $P_3$ defined in  theorem
  \ref{explicit}.
\end{itemize}
This completes the proof because we have provided a generating set of the
permutations. \fin

So  corollary \ref{bir} is proved and it implies theorem \ref{mainbir}
  from remark \ref{Hrational}. \fin

\subsection*{A normal form for   $5$ general lines in $\PP_5$}
Explicit forms of  definition \ref{defexplicit} and  theorem
\ref{explicit} have the the following straighforward  translation, that should help to
handle $5$ lines in $\PP_5$ or to understand $(G(2,V_6))^5/PGL(V_6)$.
% \begin{corollary}
%   Five lines in general position in $\PP_5$ can be put in the following form:
%   $$\epsilon_0 \wedge \epsilon_3,\  \epsilon_1 \wedge \epsilon_4,\ \epsilon_2 \wedge \epsilon_5,\  (\epsilon_0+\epsilon_1+\epsilon_2)\wedge(\epsilon_3+\epsilon_4+\epsilon_5)$$
%   $$(-\epsilon_0+v\epsilon_4-\epsilon_5)  \wedge(u.\epsilon_1-\epsilon_2+a_{1,2,4}\epsilon_3+   (a_{0,2,4}+a_{1,2,4}-a_{2,3,4})\epsilon_4+(a_{1,2,4}-a_{1,3,4}+1)\epsilon_5)$$    
% for some basis $(\epsilon_i)_{0\leq i \leq 5}$ of $V_6$, where $(u,v)$ are the roots of the
% equation: $X^{2}+X\cdot (1+a_{0,2,4}-a_{0,3,4})+a_{0,2,4}=0$. 
% \end{corollary}
\begin{corollary}\label{unirationalisation}
  Five lines in general position in $\PP_5$ can be put in the following form:
  $$\eps_0 \wedge \eps_3,\  \eps_1 \wedge \eps_4,\ \eps_2 \wedge \eps_5,\  (\eps_0+\eps_1+\eps_2)\wedge(\eps_3+\eps_4+\eps_5)$$
  $$(-\eps_0+v\eps_4-\eps_5)  \wedge(u.\eps_1-\eps_2+a_{1,2,4}\eps_3+   e_1.\eps_4+e_2.\eps_5)$$    
  for some basis $(\eps_i)_{0\leq i \leq 5}$ of $V_6$, and some complex parameters
  $u,v,a_{1,2,4},e_1,e_2$. 
\end{corollary}
\proof Let's use again  notations of  proposition \ref{P4-5sec}. From
five general lines in $\PP_5$, we obtain a five dimensional subspace $W_5\dual$
of $\bigwedge\limits^2 V_6$ containing the corresponding decomposable
elements. So choose a general four dimensional vector subspace $W_4\dual$ of
$W_5\dual$, then $(W_4\dual,W_5\dual)$ is a general element of the incidence variety
$I_{4,5}$. So from  Theorem \ref{explicit}  and  Corollary \ref{bir} the
corresponding element of $W_5\otimes \bigwedge\limits^2 V_6$ can be written with
 notation of  definition \ref{defexplicit}: $M_{0123}+w_4.M_4$. So we
obtain the proposition. \fin

% e_1=a_{0,2,4}+a_{1,2,4}-a_{2,3,4},e_2=1+a_{1,2,4}-a_{1,3,4}
%where $(u,v)$ are the roots of the equation: $X^{2}+X\cdot (1+a_{0,2,4}-a_{0,3,4})+a_{0,2,4}=0$. 
\subsection{Questions on the magic square}
\begin{remark}
Let $X$ be a non degenerate subvariety of $\PP_{n-1}$. Then the projection of
$X$ from a general linear space of dimension $d-2$ is expected to have a finite
number $n_{d,X}$ of points of multiplicity $d$ when:
$$ d^2+d(\dim(X)-n-1)+n=0.$$
\end{remark}
Varieties related to the magic square are famous solutions of this problem for $d=2$
or $d=3$ with $n_{d,X}=1$. For these varieties, what is the number
$n_{\frac{n}{d},X}$?

For the Veronese surface we have $n_{2,X}\neq n_{3,X}$, but
for $\PP_2 \times \PP_2$, $v_3(\PP_1)$, $\PP_1 \times \PP_1 \times \PP_1$ (Cf
\cite{H}) we have $n_{d,X}=n_{\frac{n}{d},X}=1$. And now, from proposition
\ref{P4-5sec}  this equality is also true for $G(2,6)$.

\section{Applications }\label{appli}
Let $V_{10}$ be a $10$-dimensional vector space over the complex numbers. In
this section, we will first explain the relationship between two known 
constructions associated to the choice of a general element of $\bigwedge\limits^3
V_{10}$. Then we will discuss how the results of the previous section should be
related to the symplectic form of the varieties constructed in \cite{D-V}.

\subsection{Peskine's example in $\PP_9$}
This example was constructed by C. Peskine to obtain a smooth non quadratically
normal variety of codimension $3$. 

Let $\PP_9$ be a $9$ dimensional projective space over the complex numbers, and
denote by $V_{10}$ the vector space $V_{10}=H^0 (\OO_{\PP_9}(1))$. Let $\alpha$
be a general element of $\bigwedge\limits^3 V_{10}$, and denote by
$\Omega^i_{\PP_9}$ the $i$-th exterior power of the cotangent sheaf of $\PP_9$. From the 
identification $\bigwedge\limits^3 V_{10} = H^0 \Omega^2_{\PP_9}(3)$, we obtain
a skew-symmetric map $M_\alpha$ from $(\Omega^1_{\PP_9})\dual(-1)$ to
$\Omega^1_{\PP_9}(2)$ and an exact sequence:
$$
\begin{CD}
  0@>>> \OO_{\PP_9}(-3)  @>>>
  (\Omega^{1}_{\PP_9})\dual(-1)@>M_{\alpha}>>\Omega^1_{\PP_9}(2)@>>> I_{Y_\alpha}(4)
  @>>> 0 
\end{CD}
$$
where $I_{Y_\alpha}$ is the ideal of the smooth variety of dimension $6$ defined by
the $8$ by $8$ pfaffians of $M_{\alpha}$. The following proposition is  directly
deduced from the previous exact sequence.
\begin{proposition}
The variety $Y_{\alpha}$ is such that $h^1(I_{Y_{\alpha}}(2))=1$ and its canonical sheaf is
$\omega_{Y_{\alpha}}=\OO_{Y_{\alpha}}(-3)$.   
\end{proposition}

\subsection{Debarre-Voisin's manifold as a parameter space}
Denote by $G(6,V_{10}\dual)$ the Grassmannian of $6$ dimensional subspaces of
$V_{10}\dual$. Let $K_6$ (resp.  $Q_4$) be the tautological subbundle
(resp. quotient bundle). For any $p\in G(6,V_{10}\dual)$, the corresponding
$5$-dimensional projective subspace of $\PP_9$ will be denoted by $\kappa_p$.

\medskip
Debarre and Voisin proved in \cite{D-V} the following:
\begin{theorem}(\cite{D-V} Th 1.1). 
Let  $\alpha$  be a general element of $\bigwedge\limits^3 V_{10} =
H^0(\bigwedge\limits^3 K_6\dual)$. The  subvariety $Z_{\alpha}$ of
$G(6,V_{10}\dual)$ defined by  
 the vanishing locus of the section $\alpha$ of $\bigwedge\limits^3 K_6\dual$ is
 an irreducible hyper-K\"ahler manifold of dimension $4$ and second betti number $23$. 
\end{theorem}
We can now remark the following relation between $Y_{\alpha}$, $Z_{\alpha}$ and
Palatini threefolds:
\begin{proposition}
 Let $p$ be a general element of $Z_{\alpha}$. The scheme defined by the
 intersection $Y_{\alpha}\cap \kappa_p$ is a Palatini threefold.
\end{proposition}
\proof
The restriction of $\Omega^1_{\PP_9}(1)$ to $\kappa_p$ is
$\Omega^1_{\kappa_p}(1)\oplus 4\OO_{\kappa_p}$. The vanishing of the restriction of
$\alpha$ to $\kappa_p$ implies that the restriction of $M_\alpha$  to $\kappa_p$
is: $\left(\begin{array}[c]{cc}
    0 & \alpha_p\\
    -\t{\alpha_p} & \beta
  \end{array}\right)$ with respect to the direct sums: $
(\Omega^1_{\kappa_p})\dual(-1)\oplus 4 \OO_{\kappa_p} \to
(\Omega^1_{\kappa_p})(2) \oplus 4\OO_{\kappa_p}(1)$. So the ideal generated by
the pfaffians of size $8$ of this map is also the ideal generated by the maximal
minors of $\alpha_p: 4 \OO_{\kappa_p} \to (\Omega^1_{\kappa_p})(2)$. In conclusion the scheme defined by the
 intersection $Y_{\alpha}\cap \kappa_p$ is a Palatini threefold as in  remark
 \ref{resolpala}. \fin

\bigskip
Moreover, the following construction globalize  definition \ref{palatini} and
the pfaffian cubic surface over $Z_{\alpha}$.
\begin{remark}\label{Zreseau}
The restriction of the bundle $\bigwedge\limits^2 K_6\dual \otimes Q_4\dual$ to
$Z_{\alpha}$ has a non trivial section. It gives   an injective map:
$$ (Q_4)_{|Z_{\alpha}} \longrightarrow (\bigwedge\limits^2 K_6\dual)_{|Z_{\alpha}} $$
\end{remark}
\proof
The section $\alpha$ of $(\bigwedge\limits^3 K_6\dual)$ gives a map from $K_6$
to $(\bigwedge\limits^2 K_6\dual)$. But the restriction of this map to
$Z_{\alpha}$ is zero, so it induces a map from the quotient $(Q_4)_{|Z_{\alpha}}$
 to $(\bigwedge\limits^2 K_6\dual)_{|Z_{\alpha}} $. The injectivity of this maps
 of  $\OO_{Z_{\alpha}}$-modules follows from
 the assumption that $\alpha$ is general. \fin

\subsection{Conjectures on the symplectic form on $Z_{\alpha}$}
\begin{remark}\label{5P1}
Let $p$ be a general element of $Z_{\alpha}$. The tangent space
$\mathcal{T}_{(Z_{\alpha},p)}$ to     $Z_{\alpha}$ at $p$  contains a canonical
set of $5$ vector spaces of dimension $2$.
\end{remark}
\proof 
Let $p$ be a general point of $Z_{\alpha}$. From remark \ref{Zreseau}, the fiber
$Q_{4,p}$ is a $4$-dimensional subspace of $(\bigwedge\limits^2
K_{6,p}\dual)$. From proposition \ref{P4-5sec}, we obtain in
$K_{6,p}\dual$,   a canonical set of five
vector subspaces $(L_i)_{0\leq i \leq 4}$ of dimension $2$ such that
$\bigoplus\limits_{0\leq i \leq 4} \bigwedge\limits^2 L_i$ contains $Q_{4,p}$. So the restriction of
the map:
\begin{equation}
  \label{tens21vers3}
m_{21}: \bigwedge\limits^2 K_{6,p}\dual \otimes K_{6,p}\dual \to \bigwedge\limits^3 K_{6,p}\dual  
\end{equation}
gives the following commutative diagram of exact sequences:
$$
\begin{CD}
  0 @>>> \mathcal{T}_{(Z_{\alpha},p)} @>>> Q_{4,p} \otimes
  K_{6,p}\dual @>>> \bigwedge\limits^3 K_{6,p}\dual @>>> 0 \\
  @. @. @VVV @VVV \\
  @. @. (\bigoplus\limits_{0\leq i \leq 4} \bigwedge\limits^2 L_i )\otimes
  K_{6,p}\dual @>>> \bigwedge\limits^3 K_{6,p}\dual @>>>0
\end{CD}
$$
where the vertical maps are injectives and the first row is the normal
sequence of $Z_{\alpha}$ in $G(6,V_{10}\dual)$ at the point $p$. Now remark that $m_{21}$
vanishes on each $\bigwedge\limits^2 L_i \otimes L_i$ because $L_i$ has
dimension $2$. So we can identify the kernel of the second row of the previous
diagram with the $10$-dimensional vector space $ \bigoplus\limits_{0\leq i \leq
  4}\bigwedge\limits^2 L_i \otimes L_i$, and we obtain an injection:
$$\mathcal{T}_{(Z_{\alpha},p)} \hookrightarrow \bigoplus\limits_{0\leq i \leq
  4}\bigwedge\limits^2 L_i \otimes L_i .$$
So in general, the kernel of each projection $ \mathcal{T}_{(Z_{\alpha},p)} \to
\bigwedge\limits^2 L_i \otimes L_i$ gives a $2$ dimensional vector subspace of $
\mathcal{T}_{(Z_{\alpha},p)}$. \fin

\bigskip
Now we can remark that five points of $G(2,\mathcal{T}_{Z_{\alpha},p})$ should
define an hyperplane $\gamma$ in $\bigwedge\limits^2 \mathcal{T}_{Z_{\alpha},p}$. Some
random examples with \cite{Macaulay2} let us expect that the ideal of these five lines $(l_i)$
in $\PP(\mathcal{T}_{Z_{\alpha},p})$ is given by the maximal minors of the map:
$$
K_{6,p} \otimes \OO_{\PP(\mathcal{T}_{Z_{\alpha},p})} \to Q_{4,p}\otimes  \OO_{\PP(\mathcal{T}_{Z_{\alpha},p})}(1)
$$
obtained from the inclusion of the tangent space to $Z_{\alpha}$ in the tangent
space to $G(6,V_{10}\dual)$. But if  the alternate form $\gamma$ was degenerated, its
kernel would give a line in $\PP(\mathcal{T}_{Z_{\alpha},p})$ intersecting each
$l_i$. But a variety defined by quartic hypersurfaces can't have a $5$-secant
line, so we can expect the following:
\begin{conjecture}
  The five vector spaces of dimension $2$ canonically defined in the remark
  \ref{5P1} are maximal isotropic subspaces for the symplectic form
  on $\mathcal{T}_{Z_{\alpha}}$  constructed by Debarre and Voisin.
\end{conjecture}

%%%%%%%%%%%%%%%%%%%%%%%%%%%%%%%%%%%%%%%%%%%%%%%%%%%%%%%%%%%%%%%%%%%%%%%%%%%%%%%%%%%

\end{document}